\documentstyle[a4,epsf,sectioneqno]{article}
\newtheorem{theorem}{Theorem}[section]
\newtheorem{lemma}[theorem]{Lemma}
\newtheorem{cor}[theorem]{Corollary}\newtheorem{defn}[theorem]{Definition}

\newtheorem{rem}[theorem]{Remark}

  \title{ On the spectrum of second order differential operators with complex coefficients
  }
\author{    B.M.Brown,  D.K.R.M$^c$Cormack \\ 
 Department of Computer Science, \\ University of Wales,  Cardiff, PO Box 916,
Cardiff CF2 3XF, U.K. 
\and W.D.Evans \\ School of Mathematics, \\ University of Wales,  Cardiff, Senghennydd Road, Cardiff, CF2 4YH,  U.K. 
\and M Plum \\ Mathematisches Institut I,\\ Universit\"at Karlsruhe,76128  Karlsruhe, Germany
\and 
 {\it  Dedicated to the memory of Professor Dr. Friedrich Goerisch}
  }
\begin{document}
\maketitle
\begin{abstract}
The main objective of this paper is to extend the pioneering work of Sims in \cite{Sims57} on second-order linear differential equations with a complex coefficient, in which he obtains an analogue of the Titchmarsh-Weyl theory and classification.
The generalisation considered exposes interesting features not visible in the special case in \cite{Sims57}.
An $m$-function is constructed (which is either unique or a point on a 
``limit-circle") and the relationship between its properties and the spectrum of underlying 
m-accretive differential operators analysed.
The paper is a contribution to the study of non-self-adjoint operators;
in general the spectral theory of such operators is rather fragmentary, and further study is being driven by important physical applications, to hydrodynamics, electro-magnetic theory and nuclear physics, for instance.
\end{abstract}

\renewcommand{\baselinestretch}{2.0}
\large \normalsize     

 \newcommand{\beq}{\begin{equation}}
\newcommand{\enq}{\end{equation}}
\newcommand{\s}{ {\cal{x}}}
\newcommand{\lek}{ \Lambda_{\eta,K}}
\newcommand{\mek}{ m_{\eta,K}}
\newcommand{\bab}{ \stackrel{\smile}{\frown}}

\section{Introduction}
In \cite{Sims57} Sims obtained an extension of the Weyl limit-point, limit-circle classification for the differential equation
\beq
M[y] =-y^{''}+qy=\lambda y,\;\;\; \lambda \in {\bf C}, \label{eq:1.1}
\enq
on an interval $[a,b)$, where $q$ is  complex-valued, and the end-points $a,b$ are respectively regular and singular. Under the assumption that ${\rm Im} q(x) \leq 0$
for all $x \in [a,b)$, Sims proved that for $ \lambda \in {\bf C_+}$, there exists at least one solution of (\ref{eq:1.1}) 
which lies in the weighted space $L^2(a,b; {\rm Im}[ \lambda-q] dx)$; such a solution lies in $L^2(a,b)$.   There are now three distinct possibilities for  $ \lambda \in {\bf C_+} $:
(I)   there is, up to constant multiples, precisely one solution of (\ref{eq:1.1}) in $L^2(a,b; {\rm Im}[ \lambda-q] dx)$ and $L^2(a,b)$, (II)
 one solution in $L^2(a,b; {\rm Im}[ \lambda-q] dx)$ but all in $L^2(a,b)$, and (III) all in $L^2(a,b; {\rm Im}[ \lambda-q] dx)$. This classification is independent
of $\lambda \in {\bf C_+}$ and, indeed, if all solutions  of (\ref{eq:1.1})
are in $L^2(a,b; {\rm Im}( \lambda-q) dx)$ or in $L^2(a,b)$, for some $\lambda$,
 it remains so  for all $\lambda \in {\bf C }$.
  At the  core of Sims' analysis is an analogue for (\ref{eq:1.1}) of the Titchmarsh-Weyl $m$-function  whose properties determine the self-adjoint realisations of
 $-\frac{ d^2}{dx^2}
+q$ in $L^2(0,\infty)$
 when $q$ is real and appropriate boundary conditions are prescribed at $a$ and $b$. Sims made a  thorough  study of the ``appropriate" boundary conditions  and  the spectral  properties of the resulting operators in the case of complex $q$. The extension of the theory for an interval $(a,b)$ where both end points are singular follows in a standard way.
\par
We have two objectives in this paper. Firstly, we construct an  analogue
of the Sims theory to the equation
\beq
-(py^{'})^{'} + q y = \lambda w y \label{eq:1.2}
\enq
where $p$ and $q$  are both complex-valued, and $w$ is a positive weight function.
This is not simply a straightforward generalisation of \cite{Sims57}, for the 
general problem exposes problems and properties of (\ref{eq:1.2}) which are
hidden in the special case 
considered by Sims; some of these features may also be seen in \cite{BK76}
where a system of the form (\ref{eq:1.2}) with $p=\omega=1$ is considered (see Remark 2.5 below). 
Secondly, once we have our analogue of the 
Titchmarsh-Weyl-Sims $m-$function, we are (like Sims) in a position to 
define
 natural  quasi m-accretive operators generated by  $-\frac{1}{w}\{\frac{ d }{dx }(p\frac{ d }{dx })
+q\}$   in $L^2(a,b;wdx)$ and to investigate their spectral properties; these, of course, depend on the analogue of the 3 cases of Sims.
Our concern, in particular, is to relate these  spectral properties to those
of the $m-$function, in a way reminiscent of that achieved  for the case of real $p,q$ by Chaudhuri and Everitt \cite{chaudhurieveritt}.
We establish the correspondence between the eigenvalues and poles of the $m-$function, but, unlike in the self-adjoint case considered in \cite{chaudhurieveritt}, there is in general a part of the spectrum which is  inaccessible from the subset of ${\bf C}$ in which the $m-$function is initially
defined and its properties determined.
However, even within this region we are able to define an $m-$function (Definition 4.10).
\par
We are grateful to the referees for comments which have helped to improve the presentation in the paper.
\section{The limit-point, limit-circle theory}
Let
\beq
M[y]= \frac{1}{w} [ - (py^{'})^{'} + q y ] \;\;\;{\rm on}\;\;[a,b)
\label{eq:2.1}
\enq
where
\newcounter{rem1}
\begin{list}%
{( \roman{rem1} )}{\usecounter{rem1}
\setlength{\rightmargin}{\leftmargin}}
 \item
$w>0$, $p \neq 0$ a.e. on $[a,b)$ and $w, 1/p \in L^1_{loc}[a,b)$;
\item
$p,q$ are complex-valued, $q \in L^1_{loc}[a.b)$ and
\beq
Q=\overline{co} \{ \frac{q(x)}{w(x)} + r  p(x):
 x \in [a,b), \; 0 < r < \infty \} \neq {\bf C}, \label{eq:2.2}
\enq
\end{list}
where $\overline{co}$ denotes the closed convex hull.
\par
 The assumptions on $w,p,q$ ensure that $a$ is a regular end-point of the equation $M[y]= \lambda w y$. We have in mind that $b$ is a singular end-point, 
i.e. at least one of $b = \infty$ or
\begin{displaymath}
\int_a^b(w + \frac{1}{ \mid p \mid } + \mid q \mid ) dx = \infty
\end{displaymath}
holds; however the case of regular $b$ is included in  the analysis.
The conditions i) and ii) will be assumed hereafter without further mention.
\par
The complement in ${\bf C}$ of the closed convex set $Q$ has one or two connected components. For $\lambda_0 \in {\bf C} \backslash Q$, denote by $K=K(\lambda_0)$ its (unique) nearest point in $Q$ and denote  by $L=L(\lambda_0)$ the tangent to $Q$ at $K$ if it exists (which it does for almost all points on the boundary of $Q$), and otherwise any line touching $Q$ at $K$.
Then if the complex plane is subjected to a translation $z \mapsto z -K$ and a rotation through an 
appropriate angle $\eta=\eta(\lambda_0) \in (-\pi, \pi]$, the image of $L$ coincides with the imaginary axis and the images of $\lambda_0$ and $Q$ lie in the new negative and non-negative half-planes respectively:
in other words, for all $x \in [a,b)$ and $r \in (0,\infty)$
\beq
Re[ \{ r p(x) + \frac{ q(x)}{w(x)}-K \} e^{i \eta} ] \geq 0 \label{eq:2.3}
\enq
and
\begin{displaymath}
Re [ (\lambda_0-K) e^{i \eta} ] <0.
\end{displaymath}
For such {\it admissible} $K, \eta$ (corresponding to some $\lambda_0 \in  {\bf C } \backslash Q $), define the half-plane
\beq
\lek := \{ \lambda \in {\bf C}: Re[(\lambda-K)e^{i \eta}] <0 \}.
\label{eq:2.4}
\enq
Note that  for all $\lambda \in \lek$
\beq
Re[(\lambda-K)e^{i \eta}]=-\delta <0
\label{eq:2.5}
\enq
where $\delta=\delta_{\eta,K}(\lambda)$ is the distance from $\lambda$ to the boundary $\partial \lek$.
Also ${\bf C}\backslash Q$ is the union of the half-planes $\lek$ over the set $S$ of admissible values of $\eta$ and $K$.
\par
We shall initially establish the analogue of the Sims-Titchmarsh-Weyl theory on the half-planes $\lek$, but subject to the condition
\beq
Re[ e^{i \eta} \cos \alpha \;\overline{ \sin \alpha}] \leq 0
\label{eq:2.6}
\enq
for some fixed $\alpha \in {\bf C}$:
the parameter $\alpha$ appears in the boundary condition at $a$ satisfied by functions in the domain of the underlying operator (see Section 4).
Denote by $S(\alpha)$ the set $\{ (\eta,K)\in S:\; (\ref{eq:2.6})\;{\rm is \;satisfied} \}$. 
We assume throughout that   
\beq
Q(\alpha) := {\bf C} \backslash \cup_{S(\alpha)} \lek = \cap_{S(\alpha)}
(  {\bf C} \backslash \lek )  \neq \emptyset.
\label{eq:2.7a}
\enq
The set $Q(\alpha)$ is clearly closed and convex, and 
  $Q(\alpha) \supseteq Q$ in general: for the important special cases $\alpha=0, \frac{\pi}{2}$, corresponding to the Dirichlet and Neumann problems, $Q(\alpha)=Q$.
In \cite{Sims57} Sims assumes that $p=w=1$ and the values of $q$ lie in ${\bf C_-}$;
 thus $\eta =\pi/2, \; K=sup_{[a,b)}[ {\rm Im }q(x)]$, are admissible  values, and  $(\eta, K) \in S(\alpha)$
if
\begin{displaymath}
{\rm - Im }[ \cos \alpha \;\overline{ \sin \alpha}]= \sinh [ 2 {\rm Im} \alpha] \leq 0,
\end{displaymath}
the assumption made by Sims. If $\alpha$ is real, then (\ref{eq:2.6}) requires $\mid \eta \mid \leq \pi/2$ if $\alpha \in [\pi/2,\pi ]$, and 
$\mid \eta \mid \geq \pi/2$ if $\alpha \in [0,\pi/2]$.
\par
We shall prove below that the spectrum of the differential operators defined in a natural way by the problems considered lie  in the set $Q(\alpha)$.
This and related results can be interpreted as implying a restriction on the range of values of boundary condition parameter $\alpha$ permitted:
if $\alpha$ satisfies (\ref{eq:2.6}) for all $\eta$ which are such that  $(\eta,K) \in S$ for some $K \in {\bf C}$, then $Q(\alpha)=Q$.
However, if $\alpha \in {\bf C}$ is given,  it is the set $Q(\alpha)$ and not $Q$, which plays the central role in general.
   \par
Let  $\theta,\phi$ be the solutions of (\ref{eq:1.2}) which satisfy
\begin{eqnarray}
\phi(a,\lambda)=\sin \alpha, & \theta(a,\lambda)= \cos \alpha \nonumber  \\
p\phi^{'}(a,\lambda)=-\cos \alpha, & p\theta^{'}(a,\lambda)= \sin \alpha \label{eq:2.6a}
\end{eqnarray}
where   $\alpha \in {\bf C}$.
On integration by parts we have, for $a \leq Y < X < b$ and $u,v \in D(M)$ defined by
\beq
D(M)= \{y : y,py^{'} \in AC_{loc}[a,b)\}, \label{eq:2.7}
\enq
that
\beq
\int_Y^X u M[v] w dx = - puv^{'} \mid_{Y}^X + \int_Y^X(pu{'}v^{'} + q uv)dx, \label{eq:2.8}
\enq
\beq
\int_Y^X (u M[v]-vM[u])w  dx = -[u,v](X)+[u,v](Y), \label{eq:2.9}
\enq
where
\beq
[u,v](x)= p(x)(u(x)v^{'}(x)-v(x)u^{'}(x)), \label{eq:2.10}
\enq
and
\begin{eqnarray}
&&\int_Y^X (u\overline{ M[v]}-\overline{v}M[u])w  dx  \nonumber \\
&&= (pu^{'}\overline{v}-\overline{p}u \overline{v}^{'})(X)-(pu^{'}\overline{v}-\overline{p}u \overline{v}^{'})(Y) +\int_Y^X [ ( \overline{p}-p)u^{'}\overline{v}^{'}+( \overline{q}-q)u \overline{v}] dx. \label{eq:2.11}
\end{eqnarray}
\par
Let $ \psi=\theta + l \phi$  satisfy
\begin{displaymath}
\psi(X) \cos \beta +(p \psi^{'})(X) \sin \beta = 0, \;\;\; \beta \in {\bf C}.
\end{displaymath}
Then
\begin{displaymath}
 l \equiv l_X(\lambda, \cot \beta)=
-\frac{ \theta(X,\lambda) \cot \beta + p(X) \theta^{'}(X,\lambda)}{ \phi(X,\lambda) \cot \beta + p(X) \phi^{'}(X,\lambda)}.
\end{displaymath}
Let
\beq
  l_X(\lambda, z):=
-\frac{ \theta(X,\lambda) z + p(X) \theta^{'}(X,\lambda)}{ \phi(X,\lambda) z + p(X) \phi^{'}(X,\lambda)},\;\;\;\;\; z \in {\bf C}. \label{eq:2.12}
\enq
This has inverse
\beq
z= z_X(\lambda, l)=
-\frac{ p(X) \phi^{'}(X,\lambda) l + p(X) \theta^{'}(X,\lambda)}{ \phi(X,\lambda) l + \theta (X,\lambda)}. \label{eq:2.13}
\enq
For   $\eta$ satisfying (\ref{eq:2.6}), the M\"obius transformation (\ref{eq:2.12}) ( note that $p( \theta \phi^{'}-\phi \theta^{'})(X)=[ \theta,\phi](X)=-1)$ is such that, for $\lambda \in \lek, \;
z \mapsto l_X(\lambda,z)$ maps the half-plane ${\rm Re}[z e^{ i \eta}] \geq 0$ 
onto a closed disc $D_X(\lambda)$ in ${\bf C}$. To see this, set $\tilde{z}=z e^{i \eta}$ and
\beq
\tilde{l}_X(\lambda, \tilde{z})= 
 -\frac{ \theta(X,\lambda) \tilde{z} + p(X) \theta^{'}(X,\lambda)e^{i \eta}}
{ \phi(X,\lambda) \tilde{z} + p(X) \phi^{'}(X,\lambda)e^{i \eta}}
=l_X(\lambda,z). \label{eq:2.14}
\enq
This has critical point
$\tilde{z}=-e^{i \eta} p(X) \phi^{'}(X,\lambda)/\phi(X,\lambda)$,
and we require this to satisfy ${\rm Re}[ \tilde{z}]<0$. We have
\begin{displaymath}
{\rm Re}[\tilde{z}]=
-{\rm Re}[ e^{i \eta}p(X) \phi^{'}(X,\lambda) \overline{\phi}(X,\lambda)/\mid \phi(X,\lambda)\mid^2 ]
\end{displaymath}
and, from (\ref{eq:2.8})
\begin{displaymath}
\int_a^X \overline{\phi}M[\phi]w dx = -p(X) \phi^{'}(X,\lambda) \overline{\phi}(X,\lambda)-
\cos \alpha \overline{ \sin \alpha }
+ \int_a^X(p \mid \phi^{'} \mid^2 + q \mid \phi \mid^2 ) dx.
\end{displaymath}
This yields
\begin{eqnarray}
 &&\mid \phi (X,\lambda )\mid^2 Re [
e^{i \eta }p(X) \phi^{'}(X,\lambda) \overline { \phi} (X,\lambda) / \mid \phi(X,\lambda )\mid^2 ] = -{\rm Re}[ e^{i \eta} \cos \alpha \overline{ \sin \alpha}]
\nonumber \\
&+& {\rm Re}[ \int_a^X e^{i \eta}\{ \frac{p}{w} \mid \phi^{'}\mid^2 + ( \frac{q}{w}
-\lambda ) \mid \phi \mid^2 \}w ]dx
\nonumber \\
&>&0 \label{eq:2.15}
\end{eqnarray}
by (\ref{eq:2.3}). 
 Thus,  when (\ref{eq:2.6}) is satisfied, $ z \mapsto l_X(\lambda,z)$ maps
${\rm Re}[z e^{i \eta}] \geq 0$ onto $D_X(\lambda)$, a closed disc with centre
\beq
\sigma_X(\lambda) = \tilde{l}_X(\lambda, e^{-i \eta} \overline { p(X) \phi^{'}
(X,\lambda)}/\overline{\phi(X,\lambda)}). \label{eq:2.17}
\enq
Furthermore $\tilde{z}=0$ is mapped onto a point on the circle $C_X(\lambda)$ bounding $D_X(\lambda)$, namely the point
\beq
\tilde{l}_X(\lambda,0)=- \theta^{'}(X,\lambda)/\phi^{'}(X,\lambda),
\label{eq:2.18}
\enq
and a calculation gives  for the radius $\rho_X(\lambda)$ of $C_X(\lambda)$
\begin{eqnarray}
\rho_X(\lambda) &=& (2 \mid {\rm Re }[ e^{i \eta} p(X) \phi^{'}(X,\lambda)  \overline{ \phi}(X,\lambda)]\mid)^{-1}   \nonumber \\
&=& \frac{1}{2} \{ -{\rm Re} [ e^{i \eta} \cos \alpha \overline { \sin \alpha}] + \int_a^X {\rm Re }[ e^{i\eta} ( p\mid \phi^{'}\mid^2 +(q-\lambda w ) \mid \phi \mid^2]dx \} ^{-1} \label{eq:2.19}
\end{eqnarray}
by (\ref{eq:2.15}).
\par
The next step is to establish that the circles $C_X(\lambda)$ are nested as $X \rightarrow b$. Set  $\psi_l=\theta+l\phi$ so that (\ref{eq:2.13})  gives
\begin{displaymath}
z=z_X(\lambda,l)=-p(X)\psi_l^{'}(X,\lambda)/\psi_l(X,\lambda).
\end{displaymath}
We have already seen that $l=l(\lambda) \in D_X(\lambda)$ if and only if
${\rm Re}[ e^{i \eta }z_X(\lambda,l)] \geq 0$, that is,
\newline
 $ {\rm Re}[ e^{i \eta} p(X) \psi^{'}_l
(X,\lambda) \overline{ \psi}_l(X,\lambda)] \leq 0$.
As in (\ref{eq:2.15}), this can be written as
\begin{displaymath}
0 \geq Re [ e^{i \eta} \{ p(a) \psi^{'}_l(a,\lambda ) \overline{\psi}_l(a,\lambda) + \int_a^X( p \mid \psi^{'}_l\mid^2 + ( q - \lambda w ) \mid \psi_l\mid^2 ) dx \} ].
\end{displaymath}
On substituting (\ref{eq:2.6a}), this gives that $l   \in D_X(\lambda)$ if 
and only if
\begin{eqnarray}
&& \int_a^X Re [ e^{i\eta}\{ p \mid \psi^{'}_l \mid^2 + ( q - \lambda w ) \mid \psi_l \mid ^2 \}] dx  \nonumber \\
&\leq& - Re [ e^{i\eta}( \sin \alpha - l \cos \alpha)( \overline { \cos \alpha} + \overline{ l \sin } \alpha )] \nonumber \\
&=:& {\cal {A}}( \alpha, \eta; l(\lambda)) \label{eq:2.20}
\end{eqnarray}
say. Note that $l \in C_X(\lambda)$ if and only if equality holds in (\ref{eq:2.20}). In view of (\ref{eq:2.3}) and (\ref{eq:2.5}), the integrand on the left-hand side  of (\ref{eq:2.20}) is positive and so $D_Y(\lambda)  \subset D_X(\lambda)$ if $X<Y$.
Hence the discs $D_X(\lambda), \; a < X< b$ are nested, and as $X \rightarrow b$ they converge to a disc $D_b(\lambda)$ or a point $m(\lambda)$: these are the {\it limit-circle} and {\it  limit-point} cases respectively. The disc $D_b(\lambda)$ and point $m(\lambda)$ depend on $\eta$ and $K$ in general, but we shall only indicate this dependence 
explicitly when necessary for clarity.
\par
Let
\beq
\psi(x,\lambda):= \theta(x,\lambda)+ m(\lambda) \phi(x,\lambda), \;\; \lambda 
\in \lek \label{eq:2.21}
\enq
where $m(\lambda)$ is either a point in $D_b(\lambda)$ in the limit-circle case, or the limit-point otherwise. The nesting property and (\ref{eq:2.20})
imply that
\beq
\int_a^b{\rm Re}[ e^{i \eta} \{ p\mid \psi^{'} \mid^2 + ( q - \lambda w ) \mid \psi \mid ^2\}] dx  \leq  {\cal {A}}(\alpha, \eta; m(\lambda)).
\label{eq:2.22}
\enq
Moreover in the limit-point case, it follows from (\ref{eq:2.19}) that
\beq
\int_a^b{\rm Re}[ e^{i \eta}  \{ p\mid \phi^{'} \mid^2 + ( q - \lambda w ) \mid \phi \mid^2 \}] dx  = \infty,
\label{eq:2.23}
\enq
whereas in the limit-circle case the left-hand side of (\ref{eq:2.23}) is finite. Also note that, by (\ref{eq:2.5}), 
a solution $y$ of (\ref{eq:1.2}) for $\lambda \in \lek$ satisfies
\beq
\int_a^b{\rm Re}[ e^{i \eta}  \{ p \mid y^{'} \mid^2 +(q-\lambda w ) \mid y \mid^2 \}]dx  < \infty
\label{eq:2.24}
\enq
if and only if
\beq
\int_a^b{\rm Re}[ e^{i \eta}  \{ p \mid y^{'} \mid^2 +(q-K w ) \mid y \mid^2 \}]dx + \int_a^b \mid y \mid^2 w dx < \infty;
\label{eq:2.25}
\enq
in particular this yields
\beq
y \in L^2(a,b;wdx). \label{eq:2.26}
\enq
In the limit-point case there is a unique solution of (\ref{eq:1.2}) for $\lambda \in \lek$ satisfying  (\ref{eq:2.25}), but it may be that all solutions satisfy (\ref{eq:2.26}). We therefore have the following analogue of Sims' result. The uniqueness referred to in the theorem is only up to  
constant multiples.
\begin{theorem}
For $\lambda \in \lek$, $(\eta,K) \in S(\alpha)$ the Weyl circles converge either to a limit-point $m(\lambda)$ or a limit-circle $C_b(\lambda)$. The following distinct cases are possible, the first two being sub-cases of the limit-point case:
\begin{itemize}
\item
Case I : there exists a unique solution of (\ref{eq:1.2}) satisfying (\ref{eq:2.25}), and this is the only solution satisfying (\ref{eq:2.26});
\item
Case II : there exists a unique solution of (\ref{eq:1.2}) satisfying (\ref{eq:2.25}),but all solutions of(\ref{eq:1.2}) satisfy   (\ref{eq:2.26});
\item
Case III: all solutions of (\ref{eq:1.2}) satisfy (\ref{eq:2.25}) and hence (\ref{eq:2.26}).
\end{itemize}
\end{theorem}
 \begin{rem}
It follows by a standard argument involving the variation of parameters formula
(c.f.\cite[Section 3 Thm. 2]{Sims57}) that the classification of (\ref{eq:1.2}) in Theorem 2.1 is independent of $\lambda$ in the following sense:
\begin{list}%
{( \roman{rem1} )}{\usecounter{rem1}
\setlength{\rightmargin}{\leftmargin}}
 \item
if all solutions of (\ref{eq:1.2}) satisfy  (\ref{eq:2.25}) for some $\lambda^{'} \in \lek$ (i.e. Case III) then all solutions of (\ref{eq:1.2}) satisfy (\ref{eq:2.25}) for all $\lambda \in {\bf C}$;
\item
if   all solutions of (\ref{eq:1.2}) satisfy (\ref{eq:2.26}) for some $ \lambda^{'} \in {\bf C}$ then all solutions of (\ref{eq:1.2}) satisfy ( \ref{eq:2.26}) for all $ \lambda \in {\bf C}$.
\end{list}
\end{rem}
\begin{rem}
 Suppose that  $p$ is real and non-negative and that for some $\eta \in [ -\frac{\pi}{2},\frac{\pi}{2}]$ and $K \in {\bf C}$,
\beq
\theta_{K,\eta}(x) = Re [ e^{i \eta}( q(x)-Kw(x))] \geq 0\;\; a.e. \; x \in(a,b).
\label{eq:2.27}
\enq
Then the condition (\ref{eq:2.25}) in the Sims characterisation of (\ref{eq:1.2}) in Theorem 2.1 for   $ \lambda \in \lek$, $(\eta,K) \in S(\alpha)$, becomes
\beq
\cos \eta \int_a^b p\mid y^{'} \mid^2dx + \int_a^b \theta_{K\eta}(x) \mid y(x) \mid^2 dx +   \int_\alpha ^b \mid y(x)\mid^2w(x) dx < \infty.
\label{eq:2.28}
\enq
\end{rem}
 In this case Remark 2.2 (i) can  be extended to the following:
\begin{list}%
{( \roman{rem1} )}{\usecounter{rem1}
\setlength{\rightmargin}{\leftmargin}}
 \item
if for some $\lambda^{'} \in {\bf C}$ all the solutions of (\ref{eq:1.2}) satisfy (\ref{eq:2.28}); then for all $\lambda \in {\bf C}$ all solutions of (\ref{eq:1.2}) satisfy (\ref{eq:2.28});
\item
if  for some $\lambda^{'} \in {\bf C}$ all the solutions of (\ref{eq:1.2}) satisfy one of
 \beq
\cos \eta \int_a^b p \mid y^{'} \mid^2 dx < \infty \label{eq:2.29}
\enq
\beq
\int_a^b \theta_{K\eta }\mid y\mid^2dx < \infty \label{eq:2.30}
\enq
then the same applies for all $\lambda \in {\bf C}$.
\end{list}
\par
The case considered by Sims in \cite{Sims57} is when $\eta = \frac{\pi}{2}, 
K=0$ in (\ref{eq:2.27}). This overlooks the interesting features present in (\ref{eq:2.28}) when $\eta \in ( -\frac{\pi}{2},\frac{\pi}{2})$, namely, that the classification in Theorem 2.1 involves a weighted Sobolev space as well as $L^2(a,b;wdx)$.
\begin{rem}
We have not been able to exclude the possibility in Cases II and III that  there exists a solution  $y$ of (\ref{eq:1.2}) for $ \lambda \in \Lambda_{\eta_1,K_1}  \cap \Lambda_{\eta_2,K_2}$ such that 
\beq
\int_a^bRe [ e^{i\eta_1} (p \mid y^{'} \mid^2 + ( q -K_1w)\mid y \mid^2 )] dx 
+   \int_a^b \mid y \mid^2 w dx < \\\infty \label{eq:2.31}
\enq
\beq
\int_a^bRe [ e^{i\eta_2 }(p \mid y^{'} \mid^2 + ( q -K_2w)\mid y \mid^2 ) ]dx 
+   \int_a^b \mid y \mid^2 w dx = \\\infty \label{eq:2.32}
\enq
for different  values of $\eta_1,\eta_2$ and $K_1,K_2$. In Case I this is  not possible by Remark 2.2. Thus, in  Cases II and III, the classification appears to depend on $K, \eta$, even under the circumstances of Remark 2.3.
\end{rem}

\begin{rem}
In \cite{BK76} a generalisation of Weyl's limit-circles theory, which includes that of
Sims, is obtained in the case of a system of the form (\ref{eq:1.2}) with $p = 
\omega = 1, \lambda = 0$ and $Im[e^{-i \eta} \; q(x)] \leq -k < 0$. The 
existence of solutions which satisfy (\ref{eq:2.25}) is established, and it is shown 
that the analogue of Case I holds when $\eta \neq \pm \frac{\pi}{2}$.
\end{rem}

\section{Properties of $m$}
Throughout the paper hearafter we shall assume  that $(\eta,K) \in S(\alpha) $. 
We denote by $\mek(\cdot)$ the function $m(\cdot)$ defined in Section 2 on
 $\lek$ whenever there is a risk of confusion. The argument in \cite[ Section 2.2]{ECT58} and
 \cite[Theorem 3]{Sims57}
remains valid in our problem to give
\begin{lemma}
In Cases I and II, $\mek$ is analytic throughout $\lek$. In Case I   the function defined by
\beq
m(\lambda)=\mek(\lambda),\;\;\lambda \in \lek \label{eq:3.1}
\enq
is well-defined on  each, of the possible two  connected  components of ${\bf C \backslash}Q (\alpha)  = \cup_{S(\alpha)} \lek$, (see (\ref{eq:2.7a})); the restriction to a connected component is   
analytic on that set.
\par
In Case III, given $m_0 \in C_b(\lambda_0), \lambda_0 \in \lek$, there exists a function $\mek$ which is analytic in $\lek$ and $\mek(\lambda_0)=m_0$,  
moreover,  a function $\mek$ can be found such that  $\mek(\lambda) \in C_b(\lambda)$ for all $\lambda \in \lek$.
\end{lemma}
{\bf Proof}   The only part not covered by the argument in \cite[Theorem 3]{Sims57}
is that pertaining to (\ref{eq:3.1}) on ${\bf C } \backslash Q(\alpha)$ in Case I.  We need only show that $ m_{\eta_1,K_1}(\lambda)= m_{\eta_2,K_2}(\lambda)$ if $ \lambda \in \Lambda_{\eta_1,K_1} \cap \Lambda_{\eta_2,K_2}$. Since in Case I, the function in
(\ref{eq:2.21}) (now denoted by $\psi_{\eta,K}(\cdot,\lambda)$ for $ \lambda \in \lek$) is the unique solution of (\ref{eq:1.2}) in $L^2(a,b;wdx)$  it follows that
\begin{displaymath}
\psi_{\eta_1,K_1}(x,\lambda)=K(\lambda) \psi_{\eta_2,K_2}(x,\lambda)
\end{displaymath}
for some $K(\lambda)$. 
On substituting the initial conditions (\ref{eq:2.6a})
we obtain $m_{\eta_1,K_1}(\lambda)=m_{\eta_2,K_2}(\lambda)$.
\par
In Case I, if ${\bf C} \backslash Q(\alpha)$ has two connected components $C_1$, $C_2$ say and $m^{(1)}, m^{(2)}$ are the $m-$functions defined on $C_1,C_2$ respectively by Lemma 3.1,
we define $m$ on ${\bf C} \backslash Q(\alpha)$  by
\begin{displaymath}
m(\lambda)= \left \{ 
\begin{array}{cc}
m^{(1)} & \lambda \in C_1, \\
m^{(2)} & \lambda \in C_2.
\end{array}
\right .
\end{displaymath}

\begin{rem}
Let $\alpha \in \{ 0,\pi  \}$ in (\ref{eq:2.20}). Then $l \in D_X(\lambda)$ implies that ${\rm Re}[ e^{i \eta} l] \geq 0$. Thus $z \mapsto l_X(\lambda,z)$
maps the half-plane ${\rm Re}[ e^{i \eta}z] \geq0$ into itself and, in particular,
$m(\cdot)$ possesses an analogue of the Nevanlinna property enjoyed by the Titchmarsh-Weyl function in the formally symmetric case. If $\alpha= \frac{\pi}{2}$, then $l \in D_X(\lambda)$ implies that ${\rm Re}[ e^{i\eta} \overline{l}] \leq 0$.
\end{rem}
\par
The argument in \cite[  Lemma 2.3]{ECT62} requires only a slight modification to give the important lemma
\begin{lemma}
Let $\lambda, \lambda^{'} \in \lek$ and $\psi(\cdot,\lambda)=\theta(\cdot,\lambda)+m(\lambda) \phi(\cdot,\lambda)$,
where $m(\lambda)$ is either the limit point or an arbitrary point in $D_b(\lambda)$ in the limit-circle case. Then
\beq
\lim_{X\rightarrow b} [ \psi(\cdot,\lambda), \psi(\cdot,\lambda^{'})](X) \equiv
\lim_{X\rightarrow b} \{p(X)[ \psi(X,\lambda)\psi^{'}(X,\lambda^{'})- \psi^{'}(X,\lambda)\psi(X,\lambda^{'})]\} =0. \label{eq:3.2}
\enq
In Case I,  (\ref{eq:3.2}) continues to hold for all $\lambda,\lambda^{'} \in {\bf C \backslash} Q(\alpha).$
\end{lemma}
{\bf Proof }
The starting point is the observation that if  Re $[z e^{i \eta}]\geq 0,$   and hence $l_X(\lambda,z)$ in (\ref{eq:2.12}) lies on the disc $D_X(\lambda)$, then with $\psi_X=\theta+l_X\phi$
\begin{displaymath}
z\psi_X(X,\lambda) + p \psi^{'}_X(X,\lambda)=0
\end{displaymath}
and similarly for $\lambda^{'}$. Then
\begin{displaymath}
[ \psi_X(\cdot, \lambda), \psi_X(\cdot, \lambda^{'})](X)=0
\end{displaymath}
and the argument proceeds as in \cite{ECT62}.
\par
Lemma 3.3 and (\ref{eq:2.8}) yield
 \begin{cor}
For all $\lambda, \lambda^{'} \in \lek$
\beq
(\lambda^{'}-\lambda) \int_a^b \psi(x,\lambda)\psi(x,\lambda^{'})w(x)dx = m(\lambda)-m(\lambda^{'}); \label{eq:3.3}
\enq
this holds for all $\lambda,\lambda^{'} \in {\bf C \backslash }Q$ in Case I.
It follows that  in Case II and III, for a fixed $\lambda^{'} \in  \Lambda_{\eta,K}$,
\beq
m(\lambda)= \frac{ m(\lambda^{'})-(\lambda-\lambda^{'})\int_a^b \theta(x,\lambda)\psi(x,\lambda^{'})w(x) dx}
{1 + (\lambda-\lambda^{'})\int_a^b \phi(x,\lambda)\psi(x,\lambda^{'})w(x)dx}
\label{eq:3.4}
\enq
defines $m(\lambda)$ as a meromorphic function in ${\bf C}$; it has a pole at $\lambda $ if and only if
\beq
1 +(\lambda-\lambda^{'})\int_a^b \phi(x,\lambda) \psi(x,\lambda^{'})w(x)=0.
\label{eq:3.5}
\enq
\end{cor}
{\bf Proof}
The identity (\ref{eq:3.3}) follows easily  from (\ref{eq:2.9}) and Lemma 3.3. In Cases II and III, $\theta(\cdot,\lambda),\phi(\cdot,\lambda) \in L^2(a,b,wdx)$, and (\ref{eq:3.4}) is derived from (\ref{eq:3.3}) on writing $\psi(\cdot,\lambda)=
\theta(\cdot,\lambda)+m(\lambda)\phi(\cdot,\lambda)$.
\begin{theorem}
Suppose that (\ref{eq:1.2}) is in Case I. Define
\begin{eqnarray}
Q_c &:= &\overline{co} \{ \frac{ q(x)}{w(x)}+r  p(x):
x \in [c,b),  \; r \in (0,\infty)\}, \label{eq:3.6} \\
Q_b &:=& \cap_{ c \in (a,b)} Q_c , \;\;\; Q_b(\alpha) = \cap_{c \in (a,b) } 
Q_c(\alpha), \label{eq:3.7}
\end{eqnarray}
where $Q_c(\alpha)$ is the set $Q(\alpha)$ defined in (\ref{eq:2.7a}) when the 
underlying interval is $[c,b)$ rather than $[a,b)$. Then $m(\lambda)$ is 
defined throughout ${\bf C \backslash} Q (\alpha)$ and has a meromorphic 
extension to ${\bf C \backslash }Q_b(\alpha)$, with poles only in 
$  Q(\alpha) \backslash Q_b (\alpha)$.
\end{theorem}
{\bf Proof}
Let $m_c(\cdot)$ denote the limit point in the problem on $[c,b)$ with $c$ now replacing $a$ in the initial conditions (\ref{eq:2.6a}); it is defined and analytic throughout  each of the possible two connected components of ${\bf C \backslash} Q_c(\alpha)$, by Lemma 3.1. Also $\psi_c := \theta _c+ m_c \phi_c $ can be
uniquely extended to $[a,b)$ with $\psi_c(x,\cdot)$ and $p\psi_c^{'}(x,\cdot)$ analytic in ${\bf C \backslash} Q_c(\alpha) $ for fixed $x$. Since we are in Case I there exists $K(\lambda)$ such that
\begin{displaymath}
\psi(x,\lambda)= K(\lambda) \psi_c (x,\lambda).
\end{displaymath}
On substituting (\ref{eq:2.6a}), we obtain
\beq
m(\lambda) = \frac{ \sin \alpha \psi_c(a,\lambda)-\cos \alpha p \psi_c^{'}(a,\lambda)}
{ \cos \alpha \psi_c(a,\lambda)+\sin \alpha p \psi_c^{'}(a,\lambda)}.
\label{eq:3.8}
\enq
This defines $m(\lambda)$ as a meromorphic function in ${\bf C \backslash} Q_c(\alpha)$
with isolated poles at the zeros of the denominator in (\ref{eq:3.8}).
In the case $b=\infty$, $Q_b$ appears in \cite[section 35]{G65}.

\section{Operator realisations of $M$}
For $\lambda \in \lek, \; (\eta,K) \in S(\alpha)$ define
\beq
G(x,y;\lambda)= \left \{
\begin{array}{cc}
-\phi(x,\lambda)\psi(y,\lambda), & a < x<y<b, \\
-\psi(x,\lambda)\phi(y,\lambda), & a<y<x<b, 
\end{array}  \right .\label{eq:4.1}
\enq
where $\phi,\psi$ are the solutions of (\ref{eq:1.2}) in (\ref{eq:2.6a}) and 
(\ref{eq:2.21}).
Recall that $m$, and hence $\psi$, depends on $(\eta,K)$ in general, but for simplicity of notation we suppress this dependency. In Case I however, Lemma 3.1 shows that $m$ is properly defined throughout ${\bf C } \backslash Q(\alpha)$. 
In Cases II and III, we know from Theorem 3.5 that $m(\cdot)$  can be continued as
a meromorphic function throughout ${\bf C}$ (but apparently still depends on $\eta$ and $K$).
For $\lambda \in \lek$ and $f \in L^2(a,b;wdx)$ define
\beq
R_\lambda f(x):= \int_a^bG(x,y;\lambda)f(y)w(y)dx. \label{eq:4.2}
\enq
 It is readily verified that $p (R_\lambda f)^{'} \in AC_{loc}[a,b)$ and from 
\begin{displaymath}
[\phi,\psi](x)=[\phi,\psi](a)=1\;\;\;(x \in (a,b))
\end{displaymath}
(see (\ref{eq:2.9}) and (\ref{eq:2.10})) that for  a.e. $ x \in (a,b)$
\beq
(M-\lambda)R_\lambda f(x)=f(x). \label{eq:4.3}
\enq
Also, for any $\lambda^{'} \in {\bf C}$
\beq
[  R_\lambda f, \phi( \cdot,\lambda^{'})] (a)=-[ \phi(\cdot,\lambda) ,\phi(\cdot,\lambda^{'})](a) \int_a^b \psi f w dx =0. \label{eq:4.4}
\enq
Moreover, if $f$ is supported away from $b$, then, by Lemma 3.3, for any $\lambda, \lambda^{'} \in \lek$,
\begin{eqnarray}
 [ R_\lambda f, \psi(\cdot,\lambda^{'})](b) &:=& 
\lim_{X \rightarrow b}[ R_\lambda f, \psi(\cdot,  \lambda^{'})](X) \nonumber \\
&=& -\lim_{X\rightarrow b} \{ [ \psi(\cdot,\lambda), \psi(\cdot,\lambda^{'})](X)
\int_a^X\phi f w dx \}  \nonumber \\
&=&0. \label{eq:4.5}
\end{eqnarray}
In Cases II and III  (\ref{eq:4.5}) holds for all $f \in L^2(a,b,;wdx)$
since then the integral on the right-hand side remains bounded as $X \rightarrow b$ and  $\lim_{X\rightarrow b} [ \psi(\cdot, \lambda), \psi(\cdot, \lambda^{'})](X)$ is zero  by (\ref{eq:3.2}).
 In Case I (\ref{eq:4.5}) continues to be true for all $\lambda, \lambda^{'}
\in {\bf C}\backslash Q(\alpha)$.
\par
Before preceding to define the realisations of $M$ which are natural to the problem,
we  need the following theorem which provides our basic tool. In the theorem $\parallel \cdot \parallel $ denotes the $L^2(a,b;wdx)$ norm.
\par

\begin{theorem}
Let $f \in L^2(a,b;wdx)$ and $\lambda \in \lek$, $(\eta, K) \in S(\alpha)$.  Then, in every case, with $\Phi \equiv R_\lambda f,$ and $\delta = dist (\lambda, \partial \lek),$
\beq
\int_a^b {\rm Re} [ e^{i \eta}  ( p \mid \Phi ^{'} \mid^2 + ( q - K w ) \mid \Phi \mid^2 )]dx  + ( {\rm Re} [ (K-\lambda )e^{i  \eta}]-\epsilon) \int_a^b \mid \Phi \mid^2 w dx \leq \frac{1}{4 \epsilon} \int_a^b \mid f \mid ^2 w dx \label{eq:4.6}
\enq 
for any $\epsilon >0$. In particular, $R_\lambda$ is bounded and
\beq
\parallel R_\lambda f \parallel \leq \frac{1}{\delta} \parallel f \parallel.
\label{eq:4.7}
\enq
\end{theorem}
{\bf Proof}
Let $f_X=\chi_{(a,X)}f$ and $\Phi_X=R_\lambda f_X$. Then, by (\ref{eq:2.8}) and (\ref{eq:4.3})
 \begin{eqnarray*}
&& \int_a^X(p \mid \Phi^{'}_X \mid^2 + (q-\lambda w )\mid \Phi_X\mid^2 ) dx=  p \overline{\Phi_X}\Phi^{'}_X \mid_a^X + \int_a^X  \overline{\Phi}_Xfw dx \\
&&= p(X)\overline{\psi(X)} \psi^{'}(X) \mid \int_a^X \phi fwdx \mid^2 -p(a) \overline{\phi}(a) \phi^{'}(a) \mid \int_a^X \psi f w dx \mid^2 +
\int_a^X \overline{\Phi_X}fwdx \\
&&= \{ \int_a^X ( p \mid \psi^{'} \mid^2 +(q-\lambda w ) \mid \psi\mid^2 ) dx +( \overline { \cos \alpha + m \sin \alpha})(\sin \alpha - m \cos \alpha )\}
\mid \int_a^X \phi f w dx \mid^2\\
&& + \overline{\sin \alpha }\cos \alpha \mid \int_a^X \psi fw dx \mid^2  
 + \int_a^X\overline{\Phi}_Xfwdx \end{eqnarray*}
from (\ref{eq:2.8}) again, and (\ref{eq:2.6a}).
 Hence, by  (\ref{eq:2.20}) and (\ref{eq:2.22}),
\begin{eqnarray*}
&\int_a^X&{\rm Re}[ e^{i\eta}( p \mid \Phi^{'}_X\mid^2 +(q-\lambda w ) \mid \Phi_X \mid^2)]dx \\
&=& \int_a^X \{ {\rm Re} [ e^{i \eta}( p \mid \psi^{'} \mid^2 + ( q-\lambda w ) \mid \psi \mid ^2 )] dx  - {\cal A}(\alpha, \eta;m(\lambda))\}\mid \int_a^X \phi fwdx \mid^2 \\
&+& {\rm Re} [ e^{i\eta} \overline{\sin \alpha}  \cos \alpha ] \mid \int_a^X  \psi fw dx \mid^2 + {\rm Re}[ e^{i\eta} \int_a^X \overline{\Phi}_Xfw]dx \\
&\leq& \int_a^X \mid \Phi_X \mid \mid f \mid w dx \leq \epsilon \int_a^X  \mid \Phi_X \mid^2 w dx + \frac{1}{4\epsilon} \int_a^X \mid f_X \mid^2w dx,
\end{eqnarray*}
whence
\begin{eqnarray*}
\int_a^b{\rm Re} [ e^{i \eta}  ( p \mid \Phi^{'}_X \mid^2 + (q-Kw) \mid \Phi_X \mid^2 ) dx ]
+ ( {\rm Re} [ e^{i \eta}(K-\lambda)]-\epsilon ) \int_a^b \mid \Phi_X \mid^2 w dx \\
\leq \frac{1}{4\epsilon}\int_a^b \mid f_X \mid^2 w dx.
\end{eqnarray*}
 As $X \rightarrow b, \; \Phi_X(x) \rightarrow \Phi(x)$  and (\ref{eq:4.6}) follows by Fatou's lemma. We also obtain from (\ref{eq:4.6}), (\ref{eq:2.3}),(\ref{eq:2.5}  and (\ref{eq:2.8}) that
\begin{displaymath}
(\delta-\epsilon)\int_a^b \mid \Phi \mid^2 w dx \leq \frac{1}{4\epsilon} \int_a^b \mid f \mid^2w dx.
\end{displaymath}
The choice $\epsilon = \frac{\delta}{2}$ yields (\ref{eq:4.7}).
\par
Theorem 4.1 enables us to establish (\ref{eq:4.5}) for all $f \in L^2(a,b;wdx)$
in Case I (and hence in all Cases).
\begin{lemma}
For $\lambda,\lambda^{'} \in \lek, \;\; (\eta, K) \in S(\alpha),$ and $f \in L^2(a,b;wdx)$
\begin{displaymath}
[R_\lambda f, \psi(\cdot,\lambda^{'}) ](b)=0.
\end{displaymath}
\end{lemma}
{\bf Proof}
Let $f_c=\chi_{[a,c]}f$, so that as $ c \rightarrow b$ we have 
\beq
f_c \rightarrow f, \;\; R_\lambda f_c \rightarrow R_\lambda f \;{\rm in }\; L^2(a,b;wdx),
\label{eq:4.8a}
\enq
\beq
[R_\lambda f_c, \psi(\cdot,\lambda^{'})](a) \rightarrow [R_\lambda f, \psi(\cdot,\lambda^{'})](a), \label{eq:4.9}
\enq
since
\begin{displaymath}
(R_\lambda f_c)(a)=-\phi(a,\lambda)\int_a^b \psi(y,\lambda)f_c(y)w dy
\rightarrow (R_\lambda f)(a),
\end{displaymath}
\begin{displaymath}
[p(R_\lambda f_c)^{'}](a) = -p\phi^{'}(a,\lambda)\int_a^b\psi(y,\lambda)f_c(y)w(y)dy \rightarrow[p(R_\lambda f)^{'}](a),
\end{displaymath}
and, by (\ref{eq:4.5}),
\beq
[R_\lambda f_c, \psi(\cdot,\lambda^{'})](b)=0.
\label{eq:4.10}
\enq
Hence, by (\ref{eq:2.9}),
\begin{eqnarray*}
[R_\lambda f, \psi(\cdot,\lambda^{'})](X) &=& [ R_\lambda (f-f_c), \psi(\cdot,\lambda^{'})](a)+[R_\lambda f_c, \psi(\cdot,\lambda^{'})](X) \\
&+& \int_a^X \{ (\lambda-\lambda^{'}) \psi(x,\lambda^{'})R_\lambda [ f-f_c](x)
+\psi(x,\lambda^{'})[f-f_c](x)\}w(x)dx \\
\rightarrow [ R_\lambda (f-f_c), \psi(\cdot,\lambda^{'})](a) & +&
\int_a^b \{ (\lambda-\lambda^{'})\psi(x,\lambda^{'}) R_\lambda [f-f_c](x)+\psi(x,\lambda^{'})[f-f_c](x) \} w(x) dx
\end{eqnarray*}
as $X \rightarrow b$, by (\ref{eq:4.10}), 
\begin{displaymath}
\rightarrow 0
\end{displaymath} 
  by (\ref{eq:4.8a}) and (\ref{eq:4.9})
\begin{rem}
 In Cases II and III, $R_\lambda$ is obviously Hilbert-Schmidt for any $ \lambda \in \lek, \; (\eta, K) \in S(\alpha)$. 
\end{rem}
In view of Theorem 4.1 and preceding remarks, it is natural to define the following operators. Let $\lambda^{'} \in \lek, \; (\eta,K) \in S(\alpha)$, be fixed and set
\begin{eqnarray}
D(\tilde{M})&:=&\{ u:u,
pu^{'} \in AC_{loc}[a,b), u,Mu \in L^2(a,b;wdx),[u,\phi(\cdot,\lambda^{'})](a)=0\;\;{\rm and}\; [ u,\psi(\cdot,\lambda^{'})](b)=0 \}, \nonumber \\
 \tilde{M}u&:=&Mu, \;\;\; u \in D(\tilde{M}). \label{eq:4.8}
\end{eqnarray}
The dependence, or otherwise,  of $D(\tilde{M})$ on $\lambda^{'}$ is made clear in 
\begin{theorem}
In Case I
\beq
D(\tilde{M})= D_1 := \{ u:u, pu^{'} \in AC_{loc}[a,b),\; u,\; Mu \in 
L^2(a,b;wdx), \; (\cos \alpha) u (a) + (\sin \alpha) p(a) u^{'}(a)=0\}.
\label{eq:4.12}
\enq
In Case II and III, $D_1$ is the direct sum
\beq
D_1=D(\tilde{M}) \stackrel{.}{+} [ \phi(\cdot,\lambda^{'})] \label{eq:4.13}
\enq
where $[\cdot ]$ indicates the linear span.
\end{theorem}
{\bf Proof}
Clearly $D( \tilde{M})  \subset D_1$: note that the boundary condition at 
$a$ in (\ref{eq:4.12}) can be written as $[u,\phi(\cdot,\lambda^{'})](a)=0$.
Let $u \in D_1$, and for $\lambda^{'} \in \lek$ set $v=R_{\lambda^{'}}[(M-\lambda^{'})u]$. Then $(M-\lambda^{'})v = (M-\lambda^{'})u$
and $[ v-u, \phi(\cdot,\lambda^{'})](a)=0$.
It follows that $v-u=K_1\phi(\cdot,\lambda^{'})$ for some constant $K_1$.
In Case I, this implies that $K=0$ since $v \in D( \tilde{M})$ and 
  $\phi(\cdot,\lambda^{'}) \not \in
L^2 (a,b;wdx)$. The decomposition (\ref{eq:4.13}) also follows since the right-hand side of (\ref{eq:4.13}) is obviously in $D_1$ in Cases II and III.
\par
In the next theorem $J$ stands for the conjugation operator $u \mapsto \overline{u}$. An operator $T$ is $J-$symmetric if $JTJ \subset T^{\ast}$ and
$J-$self-adjoint if $JTJ=T^{\ast}$ (see \cite[section III.5)]{EE87}.
Also $T$ is m-accretive if Re $\lambda < 0$ implies that $\lambda \in \rho(T)$, the resolvent set of $T$,
 and $\parallel ( T-\lambda I)^{-1} \parallel  \leq \mid {\rm Re } \lambda \mid ^{-1}$.
If for some $K \in {\bf C}$ and $\eta \in (-\pi,\pi)$, $ e^{i \eta }(T-K)$ is m-accretive, we shall say that $T$ is quasi-m-accretive;
note this is slightly different to the standard notion which does not involve the rotation $e^{i\eta}$ ( cf. \cite[section III.]{EE87}).
\par
Let $\sigma(\tilde{M})$ denote the spectrum of $\tilde{M}$.
We define the essential spectrum, $\sigma_e ( \tilde{M})$, of $\tilde{M}$
to be the complement in ${\bf C}$ of the set
\begin{displaymath}
\Delta ( \tilde{M})=\{ \lambda : ( \tilde{M}-\lambda I )
 {\rm \;is \;a\; Fredholm \; operator\; and \; ind }( \tilde{M}-\lambda I)=0\}.
\end{displaymath}
Recall that a Fredholm operator $A$ is one with closed range, finite nullity nul $A$ and  finite deficiency def $A$, and ind $A$= nul $A-$def $A$.
Thus any $\lambda \in \sigma ( \tilde{M}) \backslash \sigma_e(\tilde{M})$
is an eigenvalue of finite (geometric) multiplicity.
\begin{theorem}
The operators defined in (\ref{eq:4.8}) for any $\lambda^{'} \in \lek, \; (\eta,K) \in S(\alpha)$ (or 
(\ref{eq:4.12}) in Case I) are
$J-$self-adjoint and   quasi-m-accretive, and $\sigma(\tilde{M}) \subseteq
{\bf C}\backslash \lek$. For any $\lambda \in \lek, ( \tilde{M}-\lambda)^{-1}=R_\lambda$.
\par
In Case I, $\sigma (\tilde{M}) \subseteq Q(\alpha)$ and $\sigma_e (\tilde{M})\subseteq Q_b(\alpha)$, where $Q_b(\alpha)$ is defined in (\ref{eq:3.7}): in $Q(\alpha)\backslash Q_b(\alpha)$, $\sigma(\tilde{M})$ consists only of eigenvalues of finite  geometric multiplicity.
\par
In Cases II and III, $R_\lambda$ is compact for any $\lambda \in \rho(\tilde{M})$ and $\sigma(\tilde{M})$ consists only of isolated eigenvalues (in ${\bf C } \backslash \lek$) having finite algebraic multiplicity.
\end{theorem}
{\bf Proof} 
From $JMJ=M^+$, the Lagrange adjoint of $M$, it follows that $M$ is $J$-symmetric. Since $(\tilde{M}-\lambda)^{-1}=R_\lambda$ and $\lek \subseteq \rho(\tilde{M})$ are
established in Theorem 4.1 and the preceding remarks, it follows that $\tilde{M}$ is quasi-m-accretive, and hence also $J$-self-adjoint by Theorem III 6.7 in \cite{EE87}.
\par
In Case I, Theorem 4.1 holds for any $\lambda \in {\bf C} \backslash Q(\alpha)$ and hence $\sigma(\tilde{M}) \subseteq Q(\alpha).$  
Also, by the ``decomposition principle" (see \cite[Theorem IX 9.3 and Remark IX 9.8]{EE87}) $\sigma_e(\tilde{M}) \subseteq Q_b(\alpha)$.
\par
The compactness of $R_\lambda$ for $\lambda \in \lek$ in Cases II and III is noted in Remark 4.3, and the rest of the theorem follows.
\begin{rem}
 The argument in \cite[  Theorem 35.29]{G65} can be used to prove that in Case I of 
Theorem 4.5, either $\sigma( \tilde{M}) \backslash Q_b(\alpha)$ consists of isolated points of finite algebraic multiplicity and with no limit-point outside $Q_b(\alpha)$ or else each point of  at least one of the (possible two) connected components of $Q(\alpha)\backslash Q_b(\alpha)$ is an eigenvalue. 
We now prove that the latter is not possible.
\end{rem}
\begin{theorem}
Let (\ref{eq:1.2}) be in Case I. Then $\sigma(\tilde{M}) \subseteq Q(\alpha),\; \sigma_e(\tilde{M}) \subseteq Q_b(\alpha)$
and in $Q(\alpha)\backslash Q_b(\alpha),\; \sigma(\tilde{M})$
consists only of isolated eigenvalues of finite algebraic  multiplicity,
these points being the poles of the meromorphic extension of $m$ defined in Theorem 3.5.
\end{theorem}
{\bf Proof }  Let $\lambda \in Q(\alpha) \backslash Q_b(\alpha)$ be such that the meromorphic
extension of $m$ in Theorem 3.5 is regular at $\lambda$, and for $ c \in (a,b)$,
let $\psi(\cdot,\lambda)=K(\lambda) \psi_c(\cdot,\lambda)$ in the notation of the proof of Theorem 3.5.
Then $\psi(\cdot,\lambda)=\theta(\cdot,\lambda) + m(\lambda) \phi(\cdot,\lambda)
\in L^2(a,b;wdx)$ and the operator $R_\lambda^c$ defined by
\begin{displaymath}
R^c_\lambda f(x):= -\psi_c(x,\lambda)\int_c^x \phi(y,\lambda)f(y)w(y)dy -\phi(x,\lambda) \int_x^b \psi_c(y,\lambda)f(y)w(y)dy
\end{displaymath}
is bounded on $L^2 (c,b;wdx)$ for $c$ sufficiently close to $b$ (so that $\lambda \not \in Q_c(\alpha))$, by Theorem 4.1 applied to $[c,b)$.
Moreover (\ref{eq:4.3}) and (\ref{eq:4.4}) are satisfied by $R_\lambda$, now defined for this $\lambda \in Q(\alpha)\backslash Q_b(\alpha)$, and hence if we can prove that $R_\lambda$ is bounded on $L^2(a,b;wdx)$, it will follow that $\lambda \in \rho(\tilde{M})$, whence the theorem in view of Remark 4.6. But, for any $f \in L^2 (a,b;wdx)$, it is readily verified that
\begin{displaymath}
\parallel R_\lambda f \parallel \leq {\rm const} \{
\parallel \phi \parallel_{(a,c)}\parallel \psi \parallel  +
\parallel R^c_\lambda \parallel  \} \parallel f\parallel. 
\end{displaymath}
 Hence $\lambda \in \rho( \tilde{M})$. In Lemma 4.12 below we shall prove that $m$ is analytic on $\rho( \tilde{M})$, hence any pole of $m$  in $Q(\alpha)\backslash Q_b(\alpha)$ lies in $\sigma(\tilde{M})$.
The theorem is therefore proved.
\begin{rem}
Suppose that Case I holds.
In the notation of \cite[section IX.1]{EE87} our essential spectrum $\sigma_e$
is $\sigma_{e4}$.
However, since the operator $\tilde{M}$ is $J-$self-adjoint, by Theorem 4.5, all the essential spectra $\sigma_{ek}(\tilde{M}),\;k=1,2,3,4$ defined in \cite[Section IX.1]{EE87} coincide, by \cite[Section IX.1.6]{EE87}.
Furthermore, for any $\alpha$,
$\tilde{M}$ is a $2$-dimensional extension of the closed minimal operator generated by $M$ on
\begin{displaymath}
D_0= \{ u : u, pu^{'} \in {\rm AC}_{loc}[a,b),u,Mu,\in L^2(a,b,wdx), u(a)=p(a)u^{'}(a) =0 \}
\end{displaymath}
(cf. \cite[Theorem III 10.13 and Lemma   IX  9.2]{EE87}).
It therefore follows from \cite[  IX.1, 4.2]{EE87}  that the essential spectrum $\sigma_e(\tilde{M})$ is independent of $\alpha$.
Thus in Theorem 4.7 $\sigma_e(\tilde{M})\subseteq Q_b$, since $Q_b(0)=Q_b$.
\end{rem}
\par
We now proceed to analyse the connections between the spectrum of $\tilde{M}$ and the singularities of extensions of the $m(\cdot)$ function as is done for the Sturm-Liouville problem in \cite{chaudhurieveritt}. An important observation for this analysis is
the following  lemma. In it $(\cdot,\cdot)$ denotes the $L^2(a,b;wdx)$ inner-product.
\begin{lemma}
For all $\lambda, \lambda^{'} \in \lek, \; (\eta, K) \in S(\alpha)$,
\beq
m(\lambda)=m(\lambda^{'})-(\lambda-\lambda^{'})\int_a^b\psi^2(x,\lambda^{'})w(x)dx -(\lambda-\lambda^{'})^2 ( R_\lambda \psi(\cdot,\lambda^{'}), 
\overline{\psi}(\cdot,\lambda^{'})),
\label{eq:4.11}
\enq
\beq
m(\lambda)=[\psi(\cdot,\lambda),\theta(\cdot,\lambda^{'})](a), \label{eq:4.12a}
\enq
and
\beq
\psi(\cdot,\lambda)=\psi(\cdot,\lambda^{'})+(\lambda-\lambda^{'})R_\lambda \psi(\cdot,\lambda^{'}). \label{eq:4.13a}
\enq
\end{lemma}
{\bf Proof } 
 The identity (\ref{eq:4.11}) is an immediate consequence of (\ref{eq:3.3})
and (\ref{eq:4.13a}), and (\ref{eq:4.12a}) follows from (\ref{eq:2.6a}) and (\ref{eq:2.21}). To prove (\ref{eq:4.13a}), set $u=\psi(\cdot, \lambda )-\psi(\cdot,\lambda^{'})$.
Then $ u \in D(\tilde{M})$ by Lemma 3.3 and since
\beq[\psi(\cdot,\lambda), \phi(\cdot,\lambda^{'})](a)-[\psi(\cdot,\lambda^{'}), \phi(\cdot,\lambda^{'})](a) =0.
\label{eq:4.17aa}
\enq
Also $(\tilde{M}-\lambda)u=(\lambda-\lambda^{'})\psi(\cdot,\lambda^{'})$.
This yields $u=(\lambda-\lambda^{'})R_\lambda\psi(\cdot,\lambda^{'})$ and (\ref{eq:4.13a}) is established. The lemma is therefore proved.
 \par
Motivated by (\ref{eq:4.12a}) and (\ref{eq:4.13a})  in Lemma 4.9, we have
\begin{defn}
For $\lambda^{'} \in \lek, \; (\eta, K) \in S(\alpha),$ and $R_\lambda= ( \tilde{M}-\lambda)^{-1}$, we define $m$ on $\rho(\tilde{M})$ by
\beq
m(\lambda)= [ \Psi(\cdot,\lambda ),\theta(\cdot,\lambda^{'})](a), \label{eq:4.14a}
\enq
where
\beq
\Psi(\cdot,\lambda) = \psi(\cdot,\lambda^{'})+(\lambda-\lambda^{'})R_\lambda\psi(\cdot,\lambda^{'}). \label{eq:4.15a}
\enq
 \end{defn}
\begin{rem}
In Cases II and III, the points   $m(\lambda^{'})$ on the limit-circle for
$ \lambda^{'} \in \lek$ seem to depend on $\eta,K$ (see Remark \ref{eq:2.4})
and hence  so does the extension to $\rho ( \tilde{M})$ in Definition 4.6. This is not so in Case I, in view of Lemma 3.1.
\end{rem}
\begin{lemma}
Let $\lambda ^{'} \in \lek, \; (\eta, K) \in S(\alpha),$
 and define $m$ by (\ref{eq:4.14a})
 on $\rho(\tilde{M})$, where 
$R_\lambda =( \tilde{M}-\lambda)^{-1}$.
Then in (\ref{eq:4.15a})
\beq
\Psi(\cdot,\lambda)= \theta (\cdot,\lambda)+m(\lambda) \phi(\cdot,\lambda).
\label{eq:4.22a}
\enq
Also (\ref{eq:3.3}) and (\ref{eq:4.11}) hold for all $\lambda \in \rho(\tilde{M})$. Hence $m$ is analytic on $\rho(\tilde{M})$, and in Cases II and III, (\ref{eq:4.14a}) and (\ref{eq:3.4}) define the same meromorphic extension of $m$, while in Case I, (\ref{eq:4.14a}) defines the same meromorphic extension to ${\bf C} \backslash Q_b(\alpha)$  as that described in Theorem 3.5.
\end{lemma}
\par
{\bf Proof}
Since
\begin{displaymath}
(M-\lambda)\Psi(\cdot,\lambda )=[ (\lambda^{'}-\lambda)+(\lambda-\lambda^{'})]
\psi(\cdot,\lambda^{'})=0
\end{displaymath}
we have that
\begin{displaymath}
\Psi(\cdot,\lambda)=A \theta(\cdot,\lambda)+B\phi(\cdot,\lambda)
\end{displaymath}
for some constants $A,B$. On using (\ref{eq:2.6a}) and (\ref{eq:4.8})
it is readily verified that
\begin{eqnarray*}
A &=& -A[\theta(\cdot,\lambda ),\phi(\cdot,\lambda^{'})](a) \\
&=& -[\Psi(\cdot,\lambda ),\phi(\cdot,\lambda^{'})](a) \\
&=& -[\psi(\cdot,\lambda^{'}),\phi(\cdot,\lambda^{'})](a)-(\lambda-\lambda^{'}) [ R_\lambda\psi(\cdot,\lambda^{'}), \phi(\cdot,\lambda^{'})](a) \\
&=& 1,
\end{eqnarray*}
and
\begin{eqnarray*}
B &=&B [\phi(\cdot,\lambda),\theta(\cdot,\lambda^{'})](a) \\
&=& [ \Psi(\cdot,\lambda),\theta(\cdot,\lambda^{'})](a) \\
&=& m(\lambda)
\end{eqnarray*}
whence (\ref{eq:4.22a}). 
 Also, from (\ref{eq:4.15a})
\begin{eqnarray*}
 &&(\lambda-\lambda^{'})^2 ( R_\lambda \psi(\cdot,\lambda^{'}), \overline{\psi}(\cdot,\lambda^{'}))+ (\lambda-\lambda^{'})\int_a^b \psi^2(x,\lambda^{'})w(x)dx \\
&=& (\lambda-\lambda^{'})\int_a^b \Psi(x,\lambda)\psi(x,\lambda^{'})w(x)dx \\
&=&-\int_a^b \{ \Psi (x,\lambda)M \psi(x,\lambda^{'})-\psi(x,\lambda^{'})M \Psi(x,\lambda) \} w dx \\
&=& [ \Psi(\cdot,\lambda), \psi(\cdot,\lambda^{'})](b)-[ \Psi(\cdot,\lambda), \psi(\cdot,\lambda^{'})](a)
\end{eqnarray*}
by (\ref{eq:2.9})
\begin{displaymath}
= - [\Psi(\cdot,\lambda),\psi(\cdot,\lambda^{'})](a)
\end{displaymath}
by (\ref{eq:4.15a}) and since $\lambda \in \rho( \tilde{M})$,
\begin{eqnarray*}
&=& -m(\lambda)-m(\lambda^{'}) [ \Psi (\cdot,\lambda), \phi(\cdot,\lambda^{'})](a) \\
&=&  m(\lambda^{'})-m(\lambda)
\end{eqnarray*}
on account of (\ref{eq:4.14a}) and again using $\lambda \in \rho(\tilde{M})$. The lemma is therefore proved.
 \par
We now define, for $\lambda \in \rho( \tilde{M})$ and $f \in L^2(a,b;wdx)$,
\beq
\tilde{G}(x,y;\lambda) = \left \{\begin{array}{ll}
-\phi(x,\lambda) \Psi(y,\lambda) & a < x<y<b, \\
-\Psi(x,\lambda) \phi(y,\lambda) & a < y<x  <b,
\end{array}
\right . \label{eq:4.17}
\enq
\beq
\tilde{R}_\lambda f(x) := \int_a^b \tilde{G}(x,y;\lambda)f(x)w(x)dy,
\label{eq:4.18}
\enq
where $\Psi$ is   defined in (\ref{eq:4.22a}) and $m$   in Definition 4.10. Thus, for $\lambda \in {\bf C} \backslash Q(\alpha)$, $(\lambda \in \lek,
\; (\eta,K) \in S(\alpha),$ in Cases II and III), we have that $R_\lambda = \tilde{R}_\lambda$. We can say more, for 
(\ref{eq:4.3}), (\ref{eq:4.4}) and (\ref{eq:4.5}) hold for $\tilde{R}_\lambda$,
whenever $m(\lambda)$ is defined, and thus $\tilde{R}_\lambda=R_\lambda$ for every $\lambda$ which is such that $\tilde{R}_\lambda$ is bounded. 
This is true for every $\lambda$ at which $m$ is regular in Cases II and III. 
From (\ref{eq:4.14a}) and Lemma 4.12 we know that in Cases II and III $\lambda$ is a pole of $m(\lambda)$ if and only if
 $\lambda$ is an eigenvalue of $\tilde{M}$;
this is also true in Case I for $\lambda \not \in Q_b(\alpha)$.
\begin{theorem}
 In Cases II and III $\lambda_0$ is a pole of $m$ of order $s$ if and only if $\lambda_0$ is an eigenvalue of $\tilde{M}$ of algebraic multiplicity $s$.
\end{theorem}
{\bf Proof}
  For any $f\in L^2(a,b;wdx)$, $R_\lambda f(x)$ has a pole of order $s$ at $\lambda_0$ with residue
\begin{displaymath}
\left \{
\frac{1}{(s-1)!} \frac{\partial^{s-1}}{\partial \lambda^{s-1}}
[ (\lambda-\lambda_0)^sm(\lambda) \int_a^b \phi(x,\lambda) \phi(y,\lambda) f(y) w(y)dy]
\right \}_{\lambda=\lambda_0}.
\end{displaymath}
 This is of the form
\beq\sum_{j=0}^{s-1} \frac{\partial^{j}}{\partial \lambda^{j}}
\phi(x,\lambda_0) c_j(\lambda_0,f) \label{eq:4.25}
\enq
where the coefficients $c_j(\lambda_0,f)$ are linear combinations of
\beq
\int_a^b \frac{\partial^{j}}{\partial \lambda^{j}} 
\phi(y,\lambda_0) f(y)w(y)dy,\;\;\;j=0,1,...,s-1.
\label{eq:4.26}
\enq
From $(M-\lambda)\phi(\cdot,\lambda)=0$, it follows that for $j=0,1,...s-1$,
\beq
(M-\lambda_0)\phi_j=j\phi_{j-1}, \label{eq:4.25d}
\enq
\beq
(M-\lambda_0)^{j+1}\phi_j=0, \label{eq:4.26d}
\enq
where
\beq
\phi_j=\frac{\partial^j}{\partial\lambda^j}\phi(\cdot,\lambda_0),\; j=0,s-1.
\label{eq:4.27d}
\enq
It follows inductively from (\ref{eq:4.25d}), on using the variation of parameters, that
\beq
\phi_j \in L^2(a,b;wdx),\;\;j=0,1,...s-1.
\label{eq:4.28d}
\enq
 Let   $\Gamma_{\lambda_0}$ be a positively oriented  small circle enclosing $\lambda_0$ but excluding the other eigenvalues of $\tilde{M}$. We have
\beq
\frac{1}{2\pi i} \int_{\Gamma_{\lambda_0}} R_\lambda d\lambda=P_{\Gamma_{\lambda_0}}
\label{eq:4.27}
\enq
where $P_{\Gamma_{\lambda_0}}$ is a bounded operator of finite rank given by
(\ref{eq:4.25}): its range is spanned by $\phi_j,\; j=0,1,...,s-1$.
 The  identity  (\ref{eq:4.26d})  readily implies  that the functions in (\ref{eq:4.27d})
are linearly independent. Thus $P_{\lambda_0}$ is of rank $s$, and $s$ is the algebraic multiplicity of $\lambda_0$.
The functions in (\ref{eq:4.27d}) span the algebraic eigenspace of $\tilde{M}$
at $\lambda_0$ and are the generalised eigenfunctions corresponding to $\lambda_0$: they satisfy
\beq
(\tilde{M}-\lambda_0)^{j+1}\phi_j\neq 0,\;\;\; (\tilde{M}-\lambda_0)^j\phi_j =0\;\;j=0,1,...,s-1; \label{eq:4.30}
\enq
see \cite[Section III.4]{K76} and \cite{jbm61}.
In Case I, we expect Theorem 4.12 to remain true for $\lambda_0 \in Q (\alpha)\backslash Q_b(\alpha)$, but we have been unable to prove (\ref{eq:4.28d}) in this case.
\section{Examples}
\subsection{The sets $Q$ and $Q(\alpha)$}
Suppose that $[a,b)=[1,\infty)$ and the coefficients are of the form
\beq
p(x)= \mid p(x) \mid e^{i \phi},\;\;\; q(x)=q_1 x^{b_1} + i q_ 2 x ^{b_2},
\;\;\;w(x) = x^{\omega} \label{eq:5.1n}
\enq
where $\phi,\;  q_1 ,\;q_2 ,\;b_1, \;b_2, w$ are real constants. 
Then $q(x)/w(x)$, $x \in [1,\infty)$, lie on the curve
\beq
C:= \{ z \in {\bf C}:\;\; z = q_1x^{b_1-\omega} + i q_2 x^{b_2 - \omega},\;\; x\in [1,\infty)\}
\enq
The determination of the sets $Q$ and $Q(\alpha)$ is a straightforward exercise. As an illustration, we consider the case $\phi \in [ -\pi/2,\pi/2],\; q_1 <0,\; q_2 \leq 0,\; b_2 >b_1 > \omega$ in the Figures 1,2,3.
 The arrows indicate addition by $r p(x),\; 0 < r <\infty$, to the point 
$q(x)/w(x)$ on $C$, and the other shading in each figure is the fill-in 
required to produce the closed convex set $Q$. We set $z_0=q_1+i q_2$, 
$\tan \theta_0$ is the gradient of the tangent to $C$ at $z_0$, and $z_1$ the point 
on $C$ where the gradient is $\tan \phi$ when $\phi \geq \theta_0$ and $z_1 = z_0$ if $\phi < 
\theta_0$.
\par
The admissible values of $\eta$ (for an appropriate $K$) and the sets $Q(\alpha)$ 
for real values of the boundary value 
parameters $\alpha \in (-\pi,\pi]$ are as follows :
(recall that $Q(\alpha)$ is defined in (\ref{eq:2.7a}), where the 
admissible values of $\eta$ must now satisfy $sin 2 \alpha \;\; cos \eta \;\; 
\leq 0$)

\noindent Figure 1 : $ 0 < \eta \leq \pi/2 - \phi < \pi/2$;

\begin{center}
$Q(\alpha) = \left\{\begin{array}{ll}
Q & \rm{if}\;\; \alpha \in [-\pi/2,0] \cup [\pi/2,\pi], \\
{\bf C} & \rm{if} \;\;  \alpha \in (-\pi,-\pi/2) \cup (0,\pi/2). \\
\end{array}
\right. $
\end{center}

\noindent Figure 2 : $ 0 < \eta \leq \pi/2 - \phi < \pi$;

\begin{center}
$Q(\alpha) = \left\{\begin{array}{ll}
Q \;\;{\rm if}\;\; \alpha \in \{-\pi/2,0,\pi/2,\pi\}, \\
Q \cup \{z : \phi < arg(z - z_0) \leq 0 \} & \rm{if} \;\;\alpha \in (-\pi/2,0) \cup (\pi/2,\pi), \\
\{z : -\pi \leq arg(z - z_0) \leq \phi \} & \rm{if} \;\; \alpha \in (-\pi,-\pi/2) 
\cup (0,\pi/2). \\
\end{array}
\right. $
\end{center}

\noindent Figure 3 : $  \eta = \pi $;

\begin{center}
$Q(\alpha) = \left\{\begin{array}{ll}
Q & \rm{if} \;\; \alpha \in [-\pi,-\pi/2] \cup [0,\pi/2], \\
{\bf C} & \rm{if} \;\; \alpha \in (-\pi/2,0) \cup (\pi/2,\pi). \\
\end{array}
\right. $
\end{center}

\begin{figure}[htbp]
\centerline{
\epsfysize=7cm
\epsffile{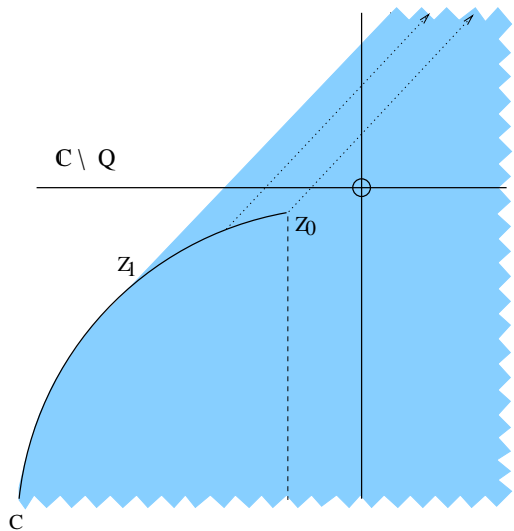}
}
\caption{$0 < \phi < \frac{\pi}{2}$}
\end{figure}

\newpage

\begin{figure}[ht]
\centerline{
\epsfysize=7cm
\epsffile{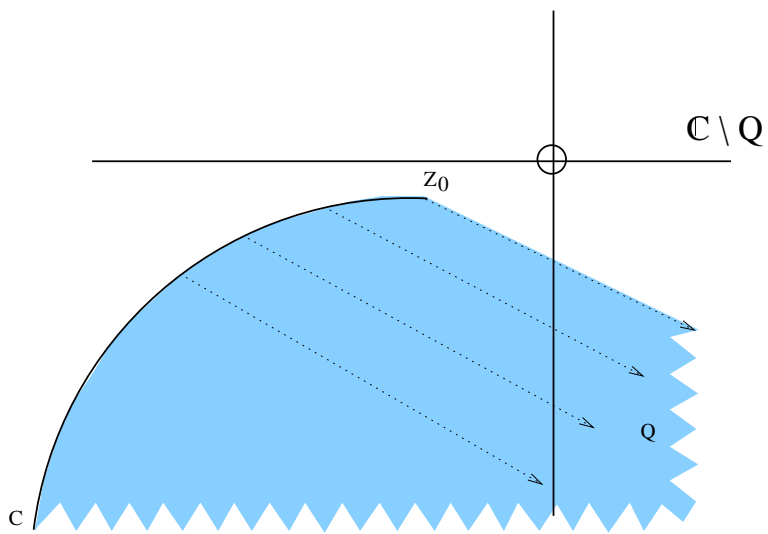}
}
\caption{$-\frac{\pi}{2} < \phi \le 0 $}
\end{figure}

\vspace{1in}

\begin{figure}[htbp]
\centerline{
\epsfysize=7cm
\epsffile{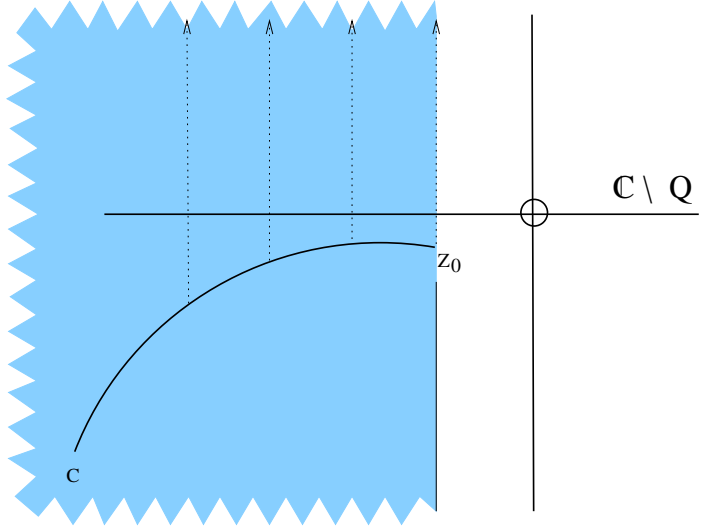}
}
\caption{$\phi = \frac{\pi}{2}$}
\end{figure}

\newpage

 \subsection{The classification of $(1.2)$}
In this section we 
analyse the Sims classification of  (\ref{eq:1.2})   when the coefficients are  
\beq
p(x)=p_1x^{ a_1}+ip_2x^{ a_2},\;q(x)=q_1x^{b_1}+iq_2x^{b_2},
\;w(x)=x^\omega,
\label{eq:5.1}
\enq 
where $p_j,q_j, a_j,b_j\;(j=1,2)$ and $\omega$ are real, and $x \in [1,\infty)$. We write $A={\rm max }\; (a_1,a_2)$ and 
$B={\rm max }\; (b_1,b_2,\omega)$. Our results follow from an   analysis of the asymptotic behaviour of
linearly independent solutions of (\ref{eq:1.2}) at infinity as given by 
the Liouville-Green formulae \cite{Eastham1}. A general description covering all 
cases is far too complicated and hardly helpful. Instead, we provide a 
prescription for determining the classification. In each specific case the details are 
straightforward, though tedious.

\subsubsection{The case $A-B <2$}

In this case, linearly independent solutions $y_\pm$ exist which are such that, 
as $x \rightarrow \infty$

\beq
  y_\pm(x)   \sim   [ p(x)s(x)]  ^{-1/4} \exp \left ( \pm \int^x _1  {\rm Re}
  [( s/p)^{1/2}]   dt \right )
\label{eq:5.2}
\enq 
\beq
  p(x)y_\pm^{'}(x)   \sim   [ p(x)s(x)]  ^{1/4} \exp \left ( \pm \int^x_1 
  {\rm Re} [( s/p)^{1/2}]  dt \right )
\label{eq:5.3}
\enq 
where
$s(x)=q(x)-\lambda w(x)$ (see \cite[page 58]{Eastham1}). We use the notation $f(x) 
\sim g(x)$ to mean that $f(x)/g(x) \rightarrow 1$ as $ x \rightarrow \infty$, 
and $f(x) \bab g(x)$ if $\mid f(x)/g(x) \mid$ is bounded above and below by positive 
constants.
Note that, for $z=re^{i \theta} \;\in {\bf C},\; 0\leq \theta < 2 \pi$, $r >0$, 
we  define the $n^{th}$ root of $z$ to be the complex number 
$r^{1/n}e^{i\theta/n}$.
\par
Suppose that for some $\lek$, $(\eta, K) \in S(\alpha)$, and 
$\lambda \in \lek$, as $x \rightarrow \infty$,
\beq
  {\rm Re} \left [ \left (\frac{ s(x)}{p(x)} \right )^{1/2}\right ]   =
D x^\tau \left ( 1 + O(\frac{1}{x^\epsilon}) \right ) \;\;\; D
\neq 0,\; \epsilon >0,
\;\; D, \tau \in {\bf R}\label{eq:5.4}
\enq
and
\beq
|p(x)s(x)| \bab x^\gamma\;\; \gamma \in {\bf R}.
\label{eq:5.5}
\enq
\par
In each of the following cases, at least one of the solutions 
$y_+$ and $y_-$ is not in $L^2(1,\infty;w dx)$, and hence (\ref{eq:1.2}) is 
in Case I :
\begin{enumerate}
\item
$\tau >-1$;
\item
$\tau =-1$ and $2 \mid D \mid + \omega -\gamma/2 + 1 \geq 0$;
\item
$\tau <-1$ and $\omega -\gamma/2 + 1 \geq 0$;.
\end{enumerate}
In all other cases when $A-B <2$, and (\ref{eq:5.4}), and (\ref{eq:5.5}) 
hold, we are either in Case II or Case III: on setting
\beq
W_\pm(x):={\rm Re} \left [
e^{i\eta} \left ( p(x) \mid y^{'}_\pm (x) \mid^2 + s(x) \mid y_\pm(x)^2 
\mid\right ) \right ]  \label{eq:5.7}
\enq
we have that Case III prevails if $W_+$ and $W_-$ are both integrable 
(which can be verified using (\ref{eq:5.2}) and (\ref{eq:5.3}) ) and Case II 
otherwise.

\subsubsection{The case $A-B=2$}
In this case  the equation  (\ref{eq:1.2}) is asymptotically of Euler type.
Here the results of \cite[page 75]{Eastham1} give, with $c=1/4( \sqrt{17}-1)$
\begin{displaymath}
\mid y_+ \;\mid \bab\; x^{ 2(A-1)c},\;\;\;
\mid py_+^{'} \mid \;\bab \;x^{ 2(A-1)(\frac{1}{2}+c)}
\end{displaymath}
and
\begin{displaymath}
\mid y_- \mid \;\bab x^{- 2(A-1)(\frac{1}{2}+c)},\;\;\;
\mid py_-^{'}\; \mid \bab x^{ -2(A-1)c}.
\end{displaymath}
At least one of the solutions $y_+$, $y_-$ is not in $L^2(1,\infty; w dx)$,
and hence (\ref{eq:1.2}) is in Case I, in each of the following cases:
\begin{enumerate}
\item
$A >1$ and $\omega + 4 ( A -1)c +1 \geq 0$;
\item
$A =1$ and $\omega   \geq -1$;
\item
$A <1$ and $\omega - 4 ( A -1)(\frac{1}{2}+c) +1 \geq 0$.
\end{enumerate}
In all other cases when $A-B =2$, we are   in   Case III when 
$W_+$ and $W_-$ defined in (\ref{eq:5.7}) are both integrable, and 
Case II otherwise.

\subsubsection{The case $A-B >2$}
Here the relevant analysis is that in \cite[page 78]{Eastham1}.
It follows that 
\begin{displaymath}
\mid y_+ \mid \;\bab 1,\;\;\; \mid y_+^{'}\mid \; \bab \; x^{ (B-A)/2},
\end{displaymath}
\begin{displaymath}
\mid y_- \mid \; \bab \; x^{-(A+B)/2},\;\;\; \mid y_-^{'}\mid \; \bab \;x^{  -A }.
\end{displaymath}
\par
At least one of the solutions $y_+$,$y_-$ is not in $L^2(1,\infty; w dx)$, 
and hence (\ref{eq:1.2}) is in Case I, if $\omega - \rm{min} \{0,A+B\} \geq -1$.
If $\omega - \rm{min} \{0,A+B\} <-1$, (\ref{eq:1.2}) is in Case III if $W_\pm$ are both integrable and Case II otherwise.
\par
The case $p=w=1$ is covered in detail in \cite[Theorem III, 10.28]{EE87};
this includes the original example of Sims \cite[p. 257]{Sims57} establishing the existence of Case II.
  
\subsection{The spectra}

Finally, we investigate the spectra of the operators $\tilde{M}$ generated
in $L^2(0,\infty)$ by expressions $M$ of the form
 \begin{equation}
M[y]=-y^{''}+c x^\beta y, \;\;\;0\leq  x<\infty, \label{eq:1}
\end{equation}
where $\beta >0$ and $c \in {\bf C}$ with arg$ \;c \in [0,\pi]$;
the case arg $ c \in (\pi,2 \pi)$ is similar.
\par
If arg $c \neq \pi$, we have
\beq
Q=\{z:0\leq {\rm arg} z \leq {\rm arg}\; c\},\;\;\; Q_\infty = \emptyset.
\label{eq:2}
\enq
Suppose that
\beq
{\rm Im } [ \overline { \sin\alpha} \cos \alpha ] \geq 0.
\label{eq:3}
\enq
Then, (\ref{eq:2.6}) is satisfied for $\eta =-\frac{\pi}{2}$ and, for any 
$K>0$, $(-\pi/2,K) \in S(\alpha)$.
Consequently
\beq
Q(\alpha) \subseteq {\bf C}\backslash \Lambda_{-\pi/2,K}= \overline{\bf C_+}
\label{eq:4}
\enq
and, similarly,
\beq
Q_\infty(\alpha)=\emptyset \label{eq:5}
\enq
(see (\ref{eq:3.7})).
Also, it follows from Section 5.2.1 (item 1) that Case I holds.
Hence, by Theorem 4.7 and Remark 4.8, for arg $c\neq \pi$,   the operator realisation $\tilde{M}$ of $M$ defined in (\ref{eq:4.12})
has empty essential spectrum $\sigma_e(\tilde{M})$.
Such a result is given in \cite[Theorem 30]{G65} for the analogous problem on $(-\infty, \infty)$.
\par
If arg $c=\pi$, we have
\beq
Q=Q_\infty={\bf R} \label{eq:6}
\enq
and, if (\ref{eq:3}) is satisfied,
$
Q(\alpha) \subseteq ({\bf C} \backslash \Lambda_{-\pi/2,K}) \cap 
({\bf C }
\backslash \Lambda_{\pi/2,K})={\bf R}$, and hence
\beq
Q(\alpha) = Q_{\infty}(\alpha) = {\bf R}. \label{eq:8}
\enq

For $\lambda = i$ and $\eta = \pm \pi/2$, we now have $\mid W_\pm \mid = 
\mid y_\pm \mid ^2$ and in Section 5.2.1

\begin{displaymath}
y_\pm(x) \bab
\left \{
\begin{array}{cc}
 x^{-\beta/4} & {\rm if\;} \beta >2, \\ 
x^{-\frac{1}{2} \mp \frac{1}{2 \mid c \mid ^{1/2}}} & {\rm if \;} \beta = 2, \\  
x^{-\frac{\beta}{4}} exp [ \mp \frac{x^{1-\beta/2}}{\mid c \mid ^{1/2} 
(2-\beta)} ] & {\rm if \;} \beta < 2. 
\end{array}
\right .
\end{displaymath}
It follows that Case I holds if $\beta \leq 2$ and Case III if $\beta > 2$;
note that Case III is now the Weyl limit-circle case since $M$ is formally 
symmetric.
Hence, if arg $c=\pi$, by Theorem 4.5,
\beq
\sigma_e(\tilde{M})
\left \{
\begin{array}{cc}
= \emptyset & {\rm if\;} \beta >2, \\
\subseteq {\bf R} &{\rm if \;} \beta \leq 2. 
\end{array}
\label{eq:9}
\right .
\enq
If $\alpha$ is real, (\ref{eq:3}) is satisfied. In this case, when $\beta 
\leq 2$, $M$ is in the Weyl limit-point 
case at $\infty$ (so that $\tilde{M}$ is self-adjoint)
and $\sigma_e(\tilde{M})={\bf R}$ (see \cite[Theorem V.5.10]{ECT62}).
\par
 In Case I the identity (\ref{eq:5.17a}) below (which holds for (\ref{eq:1.2})
 in general) is often useful and reinforces Remark 4.8.
Denote the functions $\theta,\; \phi$ in (\ref{eq:2.6a})
by $\theta_\alpha, \; \phi_\alpha$ respectively,
and the corresponding $m-$function by $m_\alpha$.
Since $\alpha=0,\pi/2$ satisfy (\ref{eq:2.6})
for any $\eta$, we have $Q(0)=Q(\pi/2)=Q$.
Also, for $\lambda \not \in Q(\alpha)$, there exist $K \neq 0$ such that
\begin{displaymath}
\theta_\alpha (x,\lambda)+m_\alpha (\lambda) \phi_\alpha(x,\lambda)
=K[ \theta_{\pi/2}(x,\lambda)+m_{\pi/2}(\lambda) \phi_{\pi/2}(x,\lambda)].
\end{displaymath}
On substituting (\ref{eq:2.6a}) we have
\beq
m_\alpha(\lambda) =
\frac{ m_{\pi/2}(\lambda) \sin \alpha-\cos \alpha}
{m_{\pi/2}(\lambda) \cos \alpha + \sin \alpha}.
\label{eq:5.17a}
\enq
Hence, if $m_{\pi/2}$ is meromorphic in ${\bf C}$,
the same is true of $m_\alpha$, for any $\alpha$.
\par
An important special case of (\ref{eq:1}) is the expression for the harmonic oscillator
\begin{displaymath}
M[y]= -y^{''} +c x^2 y \;\;\;\; 0 \leq x < \infty.
\end{displaymath}
On setting
$x=\frac{z }{\sqrt{2} c^{1/4}}$, the equation $(M-\lambda)[y]=0$ becomes
\beq
-y^{''} +\frac{1}{4}z^2 y = \mu y \label{eq:10}
\enq
 where $'$ now denotes differentiation with respect to   $z$   along the ray  with argument $\frac{1}{4} {\rm arg }\; c$, and $\mu= \frac{\lambda}{2 \sqrt{c}}$.
From \cite[page 341]{WW15},  for $0 \leq {\rm arg } c < \pi$,
the unique solution of (\ref{eq:10}) in $L^2(0,\infty)$ is the parabolic
cylinder function 
  $D_{\mu-1/2}(z)$.
It follows from (\ref{eq:2.6a}) and the fact that our function $\psi$ in 
(\ref{eq:2.21}) must be a 
constant multiple of $D_{\mu-1/2}(z)$ that
\begin{displaymath}
m_{\pi/2}(\lambda) = \frac{ D_{\mu-1/2}(0)}{D^{'}_{\mu-1/2}(0)}
\end{displaymath}
and this gives
 \begin{equation}
m_{\pi/2}(\lambda) = -\frac{1}{2 c^{1/4}} \frac{\Gamma(1/4 - \frac{\lambda}{ 4 \sqrt{c}}
)}{ \Gamma(3/4 - \frac{\lambda}{ 4 \sqrt{c}})}.  
\end{equation}
This is meromorphic with poles at
 \begin{displaymath}
\lambda_n = ( 4 n+1)\sqrt{c},\;\;\; n=0,1,2,...
\end{displaymath}
When arg $c=\pi$, $Q=Q_\infty ={\bf R}$, and for $\alpha=\pi/2$,
there are  $m$-functions defined in ${\bf C_+}$ and ${\bf C_- }$ :
\begin{eqnarray*}
m^{(1)}_{\pi/2} &=&   -\frac{e^{- i\pi/4}}{2 \mid c\mid ^{1/4}} \frac{\Gamma(1/4 - \frac{\lambda}{ 4 \sqrt{c}}
)}{ \Gamma(3/4 - \frac{\lambda}{ 4 \sqrt{c}})} \;\;\; (\lambda \in {\bf  C_+}) \\
 m^{(2)}_{\pi/2}
 &=&   -\frac{e^{i\pi/4}}{2 \mid c\mid ^{1/4}} \frac{\Gamma(1/4 +\frac{\lambda}{ 4 \sqrt{c}}
)}{ \Gamma(3/4 + \frac{\lambda}{ 4 \sqrt{c}})} \;\;\; (\lambda \in  {\bf C_-}); 
\end{eqnarray*}
${\bf  C_+},{\bf C_-}$ are the connected components $C_1,C_2$ referred  to in Lemma 3.1 and the following comment.
These functions are not analytic continuations of each other and the self adjoint operator $\tilde{M}$, with $\alpha = \pi/2$, has $\sigma_e(\tilde{M})={\bf R}$: this is therefore true for all values of $\alpha$ by Remark 4.8.
Criteria on $q$ for $\sigma_e(\tilde{M}) \supseteq [0,\infty)$ in the case $p=w=1$ are given in \cite{K79}; see also \cite[Chapter VII]{G65}.

   \bibliographystyle{plain}
\bibliography{../bib_dir/bib1,../bib_dir/bibliography,../bib_dir/help,../bib_dir/specon,../bib_dir/specon2,../bib_dir/complex}

  \end{document}